\documentclass[12pt]{article}
\usepackage{amssymb}
\usepackage{amsfonts}
\usepackage{latexsym}
\usepackage{amsthm}

\newcommand{\er}[1]{{\rm(\ref{#1})}}
\def\lb{\label}

\theoremstyle{plain}
\newtheorem{theorem}{\bf Theorem}[section]
\newtheorem{lemma}[theorem]{\bf Lemma}

\newtheorem{proposition}[theorem]{\bf Proposition}
\theoremstyle{remark}

\begin{document}

\def\a{\alpha}  \def\cA{{\cal A}}     \def\bA{{\bf A}}  \def\mA{{\mathscr A}}
\def\b{\beta}   \def\cB{{\cal B}}     \def\bB{{\bf B}}  \def\mB{{\mathscr B}}
\def\g{\gamma}  \def\cC{{\cal C}}     \def\bC{{\bf C}}  \def\mC{{\mathscr C}}
\def\G{\Gamma}  \def\cD{{\cal D}}     \def\bD{{\bf D}}  \def\mD{{\mathscr D}}
\def\d{\delta}  \def\cE{{\cal E}}     \def\bE{{\bf E}}  \def\mE{{\mathscr E}}
\def\D{\Delta}  \def\cF{{\cal F}}     \def\bF{{\bf F}}  \def\mF{{\mathscr F}}
\def\c{\chi}    \def\cG{{\cal G}}     \def\bG{{\bf G}}  \def\mG{{\mathscr G}}
\def\z{\zeta}   \def\cH{{\cal H}}     \def\bH{{\bf H}}  \def\mH{{\mathscr H}}
\def\e{\eta}    \def\cI{{\cal I}}     \def\bI{{\bf I}}  \def\mI{{\mathscr I}}
\def\p{\psi}    \def\cJ{{\cal J}}     \def\bJ{{\bf J}}  \def\mJ{{\mathscr J}}
\def\vT{\Theta} \def\cK{{\cal K}}     \def\bK{{\bf K}}  \def\mK{{\mathscr K}}
\def\k{\kappa}  \def\cL{{\cal L}}     \def\bL{{\bf L}}  \def\mL{{\mathscr L}}
\def\l{\lambda} \def\cM{{\cal M}}     \def\bM{{\bf M}}  \def\mM{{\mathscr M}}
\def\L{\Lambda} \def\cN{{\cal N}}     \def\bN{{\bf N}}  \def\mN{{\mathscr N}}
\def\m{\mu}     \def\cO{{\cal O}}     \def\bO{{\bf O}}  \def\mO{{\mathscr O}}
\def\n{\nu}     \def\cP{{\cal P}}     \def\bP{{\bf P}}  \def\mP{{\mathscr P}}
\def\r{\rho}    \def\cQ{{\cal Q}}     \def\bQ{{\bf Q}}  \def\mQ{{\mathscr Q}}
\def\s{\sigma}  \def\cR{{\cal R}}     \def\bR{{\bf R}}  \def\mR{{\mathscr R}}
\def\S{\Sigma}  \def\cS{{\cal S}}     \def\bS{{\bf S}}  \def\mS{{\mathscr S}}
\def\t{\tau}    \def\cT{{\cal T}}     \def\bT{{\bf T}}  \def\mT{{\mathscr T}}
\def\f{\phi}    \def\cU{{\cal U}}     \def\bU{{\bf U}}  \def\mU{{\mathscr U}}
\def\F{\Phi}    \def\cV{{\cal V}}     \def\bV{{\bf V}}  \def\mV{{\mathscr V}}
\def\P{\Psi}    \def\cW{{\cal W}}     \def\bW{{\bf W}}  \def\mW{{\mathscr W}}
\def\o{\omega}  \def\cX{{\cal X}}     \def\bX{{\bf X}}  \def\mX{{\mathscr X}}
\def\x{\xi}     \def\cY{{\cal Y}}     \def\bY{{\bf Y}}  \def\mY{{\mathscr Y}}
\def\X{\Xi}     \def\cZ{{\cal Z}}     \def\bZ{{\bf Z}}  \def\mZ{{\mathscr Z}}
\def\O{\Omega}
\def\ve{\varepsilon}
\def\vt{\vartheta}
\def\vp{\varphi}

\def\Z{{\Bbb Z}}
\def\R{{\Bbb R}}
\def\C{{\Bbb C}}
\def\T{{\Bbb T}}
\def\N{{\Bbb N}}
\def\dD{{\Bbb D}}

\def\sign{\mathop{\rm sign}\nolimits}
\def\sign{\mathop{\rm sign}\nolimits}
\def\dist{\mathop{\rm dist}\nolimits}
\def\Ker{\mathop{\rm Ker}\nolimits}
\def\bs{\backslash}
\def\iint{\int\!\!\!\int}
\def\bl{\biggl}   \def\br{\biggr}

\def\ma{\left(\begin{array}{cc}}
\def\am{\end{array}\right)}
\def\intl{\int\limits}
\def\iintl{\iint\limits}
\def\tminf{{\underset{m\to\infty}\to \longrightarrow}}

\def\qqq{\qquad}
\def\qq{\quad}
\let\ge\geqslant
\let\le\leqslant
\let\geq\geqslant
\let\leq\leqslant

\def\ma{\left(\begin{array}{cc}}    \def\am{\end{array}\right)}
\def\iint{\int\!\!\!\int}
\def\lt{\biggl}                     \def\rt{\biggr}
\let\ge\geqslant                   \let\le\leqslant
\def\[{\begin{equation}}            \def\]{\end{equation}}
\def\wt{\widetilde}                 \def\pa{\partial}
\def\sm{\setminus}                  \def\es{\emptyset}
\def\no{\noindent}                  \def\ol{\overline}
\def\iy{\infty}                     \def\ev{\equiv}
\def\/{\over}
\def\ts{\times}
\def\os{\oplus}
\def\ss{\subset}
\def\h{\hat}
\def\wh{\widehat}
\def\Ra{\Rightarrow}
\def\ra{\rightarrow}
\def\la{\leftarrow}
\def\da{\downarrow}
\def\ua{\uparrow}
\def\lra{\leftrightarrow}
\def\Lra{\Leftrightarrow}
\def\Re{\mathop{\rm Re}\nolimits}
\def\Im{\mathop{\rm Im}\nolimits}
\def\supp{\mathop{\rm supp}\nolimits}
\def\sign{\mathop{\rm sign}\nolimits}
\def\Ran{\mathop{\rm Ran}\nolimits}
\def\Ker{\mathop{\rm Ker}\nolimits}
\def\Tr{\mathop{\rm Tr}\nolimits}
\def\const{\mathop{\rm const}\nolimits}
\def\Wr{\mathop{\rm Wr}\nolimits}

\def\th{\theta}
\def\dlint{\displaystyle\int\limits}
\def\iintt{\mathop{\int\!\!\int\!\!\dots\!\!\int}\limits}
\def\intt{\mathop{\int\int}\limits}
\def\lim{\mathop{\rm lim}\limits}
\def\mult{\!\cdot\!}
\def\BBox{\hspace{1mm}\vrule height6pt width5.5pt depth0pt \hspace{6pt}}
\def\1{1\!\!1}
\newcommand{\bwt}[1]{{\mathop{#1}\limits^{{}_{\,\bf{\sim}}}}\vphantom{#1}}
\newcommand{\bhat}[1]{{\mathop{#1}\limits^{{}_{\,\bf{\wedge}}}}\vphantom{#1}}
\newcommand{\bcheck}[1]{{\mathop{#1}\limits^{{}_{\,\bf{\vee}}}}\vphantom{#1}}
\def\nh{\bhat}
\def\nc{\bcheck}
\newcommand{\oo}[1]{{\mathop{#1}\limits^{\,\circ}}\vphantom{#1}}
\newcommand{\po}[1]{{\mathop{#1}\limits^{\phantom{\circ}}}\vphantom{#1}}
\def\ctg{\mathop{\rm ctg}\nolimits}
\def\notto{\to\!\!\!\!\!\!\!/\,\,\,}

\def\pgbrk{\pagebreak}

\title {A priori estimates for conformal mappings
on complex plane with parallel slits}

\author{ Pavel Kargaev
\begin{footnote}
{Faculty of Math. and Mech. St-Petersburg State University}
\end{footnote}
and Evgeny Korotyaev
\begin{footnote}
{Institut f\"ur  Mathematik,  Humboldt Universit\"at zu Berlin,
Rudower Chaussee 25, 12489, Berlin, Germany,
e-mail: evgeny@math.hu-berlin.de \ \
To whom correspondence should be addressed}
\end{footnote}
}
\maketitle

\begin{abstract}
\no We study the properties of a conformal mapping $z(k)$ from the
plane without vertical slits $\G_n=[u_n-ih_n, u_n+ih_n], n\in\Z$
and $h=(h_n)_{n\in\Z}\in \ell^2$, onto the complex plane without
horizontal slits $\g_n\ss\R, n\in\Z$, with the asymptotics
$z(iv)=iv+ o(1),\ v\to\iy$. Here $u_{n+1}-u_n\ge 1, n\in
\Z$. Introduce the sequences $l=(|\g_n|)_{n\in\Z}$.
We obtain a priori two-sided estimates for $\|h\|_{p,\o}, \|l\|_{p,\o}$,
where 
%
the norm $\|h\|_{p,\o}^p=\sum \o_n|h_n|^p, 1\le p\le 2$
with any weight $\o_n\ge 1, n\in \Z$.
Moreover, we determine other estimates.
\end{abstract}


\section{Introduction and  main results}
\setcounter{equation}{0}

Consider a conformal mapping $z:K_+(h)\to\C_+$ with asymptotics
$z(iv)=iv(1+o(1))$ as $v\to\iy$, where  $z(k), k=u+iv\in K(h)$.
Here the domain $K_+(h)=\C_+\cap K(h)$ for some sequence $h=(h_n)_{n\in\Z}\in\ell^{\iy}, h_n\ge 0$ and the domain $K(h)$ is given by
\[
\lb{1.1}
K(h)=\C \sm \cup_{n\in\Z}\G_n,\qq \G_n=[u_n -ih_n,u_n+ih_n],
\qq u_*=\inf_n (u_{n+1}-u_n)\ge 0,
\]
where $u_n, n\in \Z$ is strongly
increasing sequence of real numbers such that $u_n\to \pm \iy $ as
$n\to \pm \iy $. 
We fix the sequence  $u_n, n\in \Z$ and consider the conformal mapping
for various $h\in \ell^\iy$.
The difference of any two such mappings equals a real constant. Thus the imaginary part $y(k)=\Im z(k)$ is unique.
We call such mapping $z(k) $ the comb mapping.
Define the inverse mapping $k(\cdot): {\C _+}\to K_+(h)$.
It is clear that $k(z), z=x+iy\in \C _+$ has the continuous extension
into $\ol {\C _+} $. We define
"gaps" $\g_n$, "bands"  $\s_n$ and the "spectrum" $\s$ of the
comb mapping by:
$$
\g_n=(z_n^-, z_n^+)=(z(u_n-0), z(u_n+0)),\qqq
\s_n=[z_{n-1}^+, z_n^-],\qqq \s=\cup_{n\in\Z }\s_n.
$$
The function $u(z)=\Re k(z) $ is strongly  increasing on each band
$\s_n$ and $u(z)=u_n$  for all $z\in [z_n^-,z_n^+],\ n\in\Z $;
the function $v(z)=\Im k(z)$ equals zero  on each  band $\s_n$ and
is  strongly convex on each gap $\g_n\ne \es$ and has the maximum  at some point
$z_n $ given by $v(z_n)=h_n$. If  the gap  is empty we set $z_n=z_n^{\pm}$. The function $z(\cdot) $ has an analytic extension (by the symmetry)
from  the domain $K_+(h)$   onto  the domain $K(h)$
and $z(\cdot): K(h)\to z(K(h))=\cZ=\C\sm\cup \ol\g_n $ is  a conformal
mapping. These and  others properties of the comb mappings
it is possible  to  find in the papers of  Levin [Le].

For any $p\ge 1$ and  the weight $\o=(\o_n)_{n\in \Z}$, where 
$\o_n\ge 1$, we introduce the real  spaces
$$
 \ell^p_{\o}=\{f=(f_n)_{n\in \Z}:\ \ \| f \|_{p,\o}<\iy \},
\ \ \ \| f \|_{p, \o }^p=\sum_{n\in\Z } \o_n f_n^p <\iy .
$$
If the weight $\o_n=(2u_n)^{2m}, m\in \R$ for all $n\in \Z$, then we will
write $\ell_m^p$ with the norm $ \| \cdot \|_{p,m}$.
If the weight $\o_n=1$ for all $n\in \Z$, then we will write
 $\ell^p_0=\ell^p$  with the norm
 $\| \cdot \|_{p}$ and $\| \cdot \|= \| \cdot \|_{2}$.
For each $h=(h_n)_{n\in \Z}$ we introduce the sequences
$$
l=(l_n)_{n\in \Z},\ \ l_n=|\g_n|, \qq
 J=(J_n)_{n\in \Z},\qq J_n=|A_n|^{1\/2}\ge 0,\qq
A_n={2\/\pi }\int _{\g_n}v(z)dz\ge 0. 
$$
For the defocussing cubic non-linear Schr\"odinger equation
(a completely integrable infinite dimensional
Hamiltonian system), $k(z)$ is a quasi-momentum and  $A_n$ is an action variables (see \cite{K6}).  We formulate the first main result
about the estimates in terms of $\|\cdot\|_p$.

\begin{theorem}\lb{T2.1}
Let  $u_*=\inf_n (u_{n+1}-u_n)>0$. Then the following estimates hold:
\[
\lb{2.2}
\|h\|_p\le 2\|l\|_p(1+\a_p\ \|l\|_p^p), \ \ 1\le p\le 2,
\ \ \ \a_p=(2^{p+2}(2+\pi )/u_*)^p/\pi ,
\]
\[
\lb{2.3}
\|h\|_{p}\le \frac{2}{\pi}C_p^2 \|l\|_q\lt(
1+\lt[\frac{2C_p}{\pi u_*}\rt]^{\frac2{p-1}}\|l\|_q^{2\/p-1}\rt),
\ \ C_p=\lt({\pi ^2\/2}\rt)^{1/p},  \ \  p\ge 2, \ \ {1\/p}+{1\/q}=1,
\]
\[
\lb{2.4}
{\|l\|_p  \/  2}\le  \|J\|_p \le
{2\/  \sqrt{\pi }}\|l\|_p(1+\a_p\|l\|_p^p)^{1/2},\ \ \
\]
\[
\lb{2.5}
{\sqrt{\pi}  \/  2}\|J\|_p\le \|h\|_p\le
4\|J\|_p(1+\a_p 2^{p}\|J\|_p^p).
\]
\end{theorem}

Estimates \er{2.2}-\er{2.5} are new for $p\in [1,2)$.
Korotyaev \cite{K1} obtained the two-sided estimates for the case $p=2$
(see Theorem \ref{T3.5}),
for example
\[
\lb{3.7}
{1 \/  2}\|l\|\le\|h\|\le\pi \|l\|\bigg(1+{2 \/  u_*^2}\|l\|^2\bigg),
\qq if \ \ u_*>0.
\]

Introduce 
the effective masses $\mu_n^{\pm}$ for the ends $z_n^-<z_n^+$ by
\[
\lb{2.20}
z(k)-z_n^{\pm} ={(k-u_n)^2 \/  2\mu_n^{\pm}}+ O((k-u_n)^3)
\qqq {\rm as} \ \ z\to z_n^{\pm}.
\]
If $|\g_n|=0$, then we set $\m_n^{\pm}=0$. Define the sequence $\m^{\pm}=(\m_n^{\pm})_{n\in \Z}$. 
We formulate the second result.

\begin{theorem}\lb{T2.2}
Let $h\in \ell_{\o}^{p}, p\in [1,2]$ and let $u_*>0$. Then the following estimates hold:
\[
\lb{2.6}
\|h\|_{\iy }\le \min \{2\pi \|\m^{\pm}\|_{\iy },\  \|J\|_{p,\o} ,\
2\pi ^{-1/p}\|l\|_{p,\o}(1+\a_p\|l\|_{p,\o}^p)^{1/q}\},
\]
\[
\lb{2.7}
\|l\|_{p,\o}\le 2\|h\|_{p,\o} \le \x^9 \|l\|_{p,\o},\ \qq
\x=\exp {(\|h\|_{\iy }/u_*)},
\]
\[
\lb{2.8}
\|l\|_{p,\o} \le  2\|J\|_{p,\o} \le \x^5 2\|l\|_{p,\o},
\]
\[
\lb{2.9}
{\sqrt{\pi} \/  2}\|J\|_{p,\o} \le \|h\|_{p,\o}\le
\x^5\sqrt{\pi \/  2}\|J\|_{p,\o} ,
\]
\[
\lb{2.10}
\|l\|_{p,\o}\le 2\|\m^{\pm}\|_{p,\o} \le \x^{18}
\|l\|_{p,\o}.
\]
\end{theorem}

Estimates \er{2.6}-\er{2.10} are new. Korotyaev 
obtained the two-sided estimates for the space $\ell^2_m, m\ge 0$
for the even case $h_{-n}=h_n, n\in \Z$ (\cite{K2}-\cite{K4}) and
for the space $\ell^2_1$ without symmetry (\cite{K1}, (\cite{K6})).
In all these estimate the factor $\x=\exp {(\|h\|_{\iy }/u_*)}$
is absent.

\begin{proposition}\lb{T-E}
i) Estimate \er{2.2} at $p=1$ is sharp.

ii) Estimate \er{3.6} is sharp.

iii) If the estimate $\|h\|\le C\|l\|(1+\|l\|^p),\ \ p>0$
is true, then  $p\ge 1$.

\end{proposition}

 Recall that for a compact subset $\O\ss\C $ the analytic capacity
is given by
\[
\lb{2.11}
\cC=\cC (\O)=\sup\lt[|f'(\iy)|: f \ is \ analytic\ in  \  \C \sm \O;
\ \ \ |f(k)|\le 1,\ k\in\C \sm \O \rt],
\]
where $f'(\iy)=\lim_{|k|\to\iy}k(f(k)-f(\iy))$. We will use the
well known Theorem (see \cite{Iv}, \cite{Po})

\no {\bf Theorem }    {\it ( Ivanov-Pommerenke ).
 Let $E\ss\R $ be compact. Then  the analytic capacity
$\cC (E)=|E|/4$, where $|E|$ is the Lebesgue measure (the length) of the set  $E$. Moreover, the Ahlfors function $f_E$ (the unique function, which gives  $\sup$ in the definition of  the analytic capacity) has the following form:}
\[
\lb{2.12}
f_E(z)=\frac{\exp{(\frac12\f_E(z))}-1}{\exp{(\frac12\f_E(z))}+1},\ \ \
\ \ \ \  \f_E(z)=\int_E\frac{dt}{z-t}; \ \ z\in\C \sm E.
\]
We will use the following simple remark: 
Let $S_1, S_2,\dots, S_N$ be disjoint continua
in the  plane $\C ;\ \ D=\C \sm\cup_{n=1}^NS_n$. Introduce the class
$\S '(D)$ of the conformal mapping $w$ from the domain $D$ onto $\C $ with
the following asymptotics: $w(k)=k+[Q(w)+o(1)]/k,\  k\to\iy$.
If $\O\ss\C $ is compact; $ D=\C \sm \O,\ g\in\S '(D)$, then  $\cC (\O)=\cC (\C \sm g(D))$. It follows immediately from the definition of the analytic capacity.

Let $\F_+\ss \ell^\iy$ be the subset of finite sequences
of non negative numbers. Then, using  the Ivanov-Pommerenke Theorem
 and the last remark we obtain
$$
\|l (h)\|_1=\cC (\G (h)),\ \ where \ \ \  h\in\F_+, \ \G (h)=\cup [u_n-ih_n,u_n+ih_n];\ \ \ \
$$

\begin{proposition} \lb{T3.6}
Let $h\in \ell^{\iy}$ and let $u_*\ge 0$. Then
\[
\lb{2.29}
\|h\|_\iy^2\le 2Q_0={1\/\pi}\int_{\R }v(x)\,dx ,
\]
\[
\lb{3.17}
\pi Q_0\le\|h\|_{p}^{}\,\|l\|_q, \ \ \  p\ge 1,
\]
\[
\lb{3.18}
I_D\le (\frac2{\pi })^{2/p}\|h\|_{p}^{2/q}\,\|l\|_p^{2/p},
 \ \ \ 1\le  p\le 2,
\]
\[
\lb{3.19}
\pi Q_0\le\|h\|_{\iy}^{}\,\|l\|_1\le\frac2{\pi}\|l\|^2_1,
\]
\[
\lb{3.20}
\|h\|_{\iy}^{}\le\frac2{\pi}\|l\|_1,\quad \|l\|_1\le2\|h\|_1.
\]
\end{proposition}

Below we  will sometimes write $\g_n(h), z(k,h),..$, instead of $\g_n, z(k),..$, when several sequences $h\in \ell^\iy$ are being dealt with. 
Define the Dirichlet integral 
$$
I_D(h)={1\/\pi}\iint_{\C }|z'(k,h)-1|^2\,du\,dv
={1\/\pi}\iint_{\C }|k'(z,h)-1|^2\,dx\,dy,\ k=u+iv,\ z=x+iy.
$$
The last identity holds since the Dirichlet integral is invariant
under the conformal mappings.  

Now we estimate the Dirichlet integral $I_D(h) $ (or $Q_0(h)={1\/\pi}\int_{\R }v(x)dx$)
for the case $u_*\ge 0$, using the following
geometric construction. For the vector $h\in \ell^\iy,  h_n\to 0$ as $n\to\iy$, we introduce the sequence $\wt h=\wt h(h)$ by: if $h=0$, then $\wt h=0$,

if $ h\ne 0$, then 
we take an integer $n_1 $ such that $\wt h_{n_1}=h_{n_1}=\max_{n\in\Z }h_n>0$; assume that we define the numbers
 $h_{n_1},h_{n_2},\dots ,h_{n_k} $,  then we take $n_{k+1} $ such that
\[
\lb{2.15}
\wt h_{n_{k+1}}=h_{n_{k+1}}=\max_{n\in B}h_n>0,
\ \  B=\{n\in\Z : |u_n-u_{n_l}|>h_{n_l}, 1\le l\le k\}.
\]
Moreover, we let $\wt h_n=0$, if $n \notin \{n_k, k\in \Z\}$.
Now we formulate the following results

\begin{theorem}\lb{T2.4}
Let $h\in \ell^\iy, h_n\to 0 $ as $|n|\to\iy;$ and
let $\wt h=\wt h(h)$. Then the following estimates hold:
\[
\lb{2.16}
\frac1{\pi^2}\|\wt h\|^2_2\le Q_0(h)=
{I_D(h)\/2}\le\frac{2\sqrt2}\pi\|\wt h\|^2_2.
\]
\end{theorem}

We give the geometry interpretation of the estimates from all
these theorems.
Define the square differential in the domain $D=K(h)\cup\{\iy\}$ on
the Riemann sphere, which is considered as the Riemann surface with
hyperbolic boundary components by:
$$
w=(k'(z,h)-1)^2dz^2=(z'(k,h)-1)^2dk^2.
$$
Then $w$ is the analytic square differential on $D$ (in particular,
analytic at any boundary point in terms of a corresponding uniformizing
parameter). In the present paper the metric $w$ is important to get the needed
estimates.
 Moreover, these estimates have the following geometry interpretation.
The invariant length $L_n$ of the cut $\g_n$ has the following form
$$
L_n=2\int_{z_n^-}^{z_n^+}|k'(x)-1|dx=2\int_{z_n^-}^{z_n^+}
\sqrt{v'(x)^2+1}dx.
$$
Then using \er{1.3} we obtain $2h_n\le L_n\le 2(h_n+l_n)\le 6h_n$.
Moreover, the invariant area $S$ of the Riemann surface $D$ has the form
$$
S=\iint_{\C}|k'(z)-1|^2dxdy=2\pi Q_0=\int_\R v(t)dt.
$$
Using Theorem 1.1-1.5, we  estimate $S$ in terms of $(h_n)_1^N$
or $(L_n)_1^N$.

Levin [Le] 
 proved the existence of the mapping for a very general case.
First two-sided estimates for $\|h\|_{2,1}$ and $Q_2={1\/\pi}\int_\R t^2v(t)dt$ were obtained in [MO2] only for the case $u_n=\pi n, n\in \Z$.
Note that these estimates are overstated since the Bernstein inequality
was used. Garnet and Trubowitz \cite{GT} also obtained some estimates,
using the different arguments.  The authors of the present paper \cite{KK1} obtained estimates
of the various parameters, the identities  \er{1.5}. First two-sided ewstimates ( very rough)
for $\|l\|_{2,1}$ and $Q_2$ were obtained in \cite{KK2}. Identities and various sharp estimates (in terms of gap lengths and effective masses) were obtained by Korotyaev \cite{K1}-\cite{K6}.

The estimates for conformal mapping were used  to study
the inverse problem for the Schr\"odinger operator with a periodic potential \cite{KK2}, [K5-7], for the periodic  weighted operator
\cite{K8} and for the periodic Zaharov-Shabat systems \cite{K6}.

Note that the comb mappings are used in various fields of mathematics.
We  enumerate the more important directions:

{\it \no 1) the conformal mapping theory,
2) the L\"owner equation and the quadratic differentials,
3) the electrostatic problems on  the plane,
4) analytic capacity,
5)  the spectral  theory of the operators  with  periodic coefficients,
6) inverse problems for the  Hill operator and the Dirac operator,
7) KDV equation and NLS equation with periodic initial value problem.}

Finally, we shall briefly describe our motivation to derive the present results. Consider the electrostatic field in  the domain 
$K(h)=\C \sm \cup_{n\in\Z }\G_n$, where 
$\G_n=[u_n -ih_n, u_n+ih_n],  n\in \Z$,
is the system of neutral conductors, for some $h\in\ell^p_\o$
and $u_{n+1}-u_n\ge 1,\ n\in\Z$.
In other words, we embed the system of neutral conductors $\G_n, n\in \Z,$ in the  external
homogeneous electrostatic field $E_0=(0, -1)\in \R^2$ on the plane.
Then on each conductor there exists the induced charge,
positive $e_n>0$ on the lower half of the conductor $\G_n$ and negative
 $(-e_n)<0$ on the upper half of the conductor $\G_n$, since their
sum equals zero. As a result we have new perturbed electrostatic field
 $\cE\in \R^2.$ It is well known that 
$\cE=\ol {iz'(k,h)}=-\nabla y(k,h), \ k=u+iv\in K(h), \ \ z=x+iy$. 
Recall that $z(k,h)$ is the conformal mapping from $K(h)$ onto the domain
 $\cZ=\C\sm \cup \g_n$.
The function $y(k,h)$ is called the potential of the electrostatic field
in $K(h)$. The density  of the charge on the conductor has the form
$\r_e(k)=|y_u'(k,h)|/4\pi, \ k\in \G_n$ (see [LS]).
Thus we obtain the induced charge $e_n$ on the upper half of
the conductor $\G_n^+=\G_n\cap\C_+$ by:
$$
e_n={1 \/  4\pi}\int _{\G_n^+}x_v'(k)dv={1 \/  4\pi}|\g_n|.
$$
Introduce the bipolar moment $d_n$ of the conductor $\G_n$
with the charge density $\r_e(k)$ by 
$d_n={1\/4\pi}\int _{\G_n}vx_v(k)dv\ge 0$.
We transform this value into the form
$$
d_n={1 \/  2\pi}\int _{\g_n}v(x)dx={A_n\/ 4}.
$$
In the paper \cite{KK3} we study inverse problems
for the charge mapping $h\to e(h)=(e_n(h))_{n\in \Z}$ and the
bipolar moment mapping  $h\to J(h)$ acting in
$\ell^p_{\o }, p\in [1,2]$. In order to solve
the inverse problems we need a priori estimates
from Theorem \ref{T2.1} and \ref{T2.2}.

We now describe the plan of the paper. 
In Section 2 we shall obtain 
some preliminaries results and 
"local basic estimates" in Theorem \ref{T3.5}.
In Section 3 we shall prove the main theorems.
Moreover, we consider some examples, which describe
our estimates.

\section{Preliminaries}
\setcounter{equation}{0}

We recall needed results.
Below we will use very often the following simple estimate 
\[
\lb{1.3}
l_n\le 2h_n,\qq {\rm all}\qq n\in\Z,
\]
(see e.g. [MO1], [KK1]).
Hence if $h\in \ell^p_{\o}$, then  $l(h)\in \ell^p_{\o} $.

For each $h\in \ell^{\iy}$
the following estimates and identities hold
\[
\lb{1.5}
{1 \/4} \|l\|^2\le 2Q_0=I_D=\sum A_n=\|J\|^2\le 
{2\/\pi} \sum_{n\in\Z }h_nl_n,
\]
\[
\lb{2.30}
\max \lt\{ {l_n^2 \/  4}, {l_nh_n \/  \pi}\rt\}\le A_n=
{2 \/  \pi} \int _{\g_n}v(x)dx\le {2l_nh_n \/  \pi},
\]
see [KK1]. These show that functional $Q_0={1\/\pi}\int_{\R }v(x)\,dx $
is bounded for $h\in \ell^2$.
Define the effective masses $\n_n$ in the plane $K(h)$ for the  end
of the slit $[u_n+ih_n,u_n-ih_n]$ by 
\[
\lb{2.21}
k(z)-(u_n+ih_n)={(z-z_n)^2 \/  2i\n_n}+O((z-z_n)^3),
\ \ \ \ \ z\to z_n.
\]
Thus we obtain $\n_n=1/|k''(z_n)|$, if $h_n>0$ and we set $\n_n=0$
if $l_n=0$.

Below we will use the Lindel\"of principle
(see [J]), which is formulated in  the form, convenient  for us (see \cite{KK1}):

{\it  Let $h,\ \wt h\in \ell^\iy;$ and let $\wt h_n\le h_n$ for all $n\in\Z $.
Then the following estimates hold:}
\[
\lb{2.23}
y(k,\wt h)\ge y(k,h),\ \ k\in K_+(h),
\]
\[
\lb{2.26}
Q_0(\wt h)\le Q_0(h)\ \  and \ if \ Q_0(\wt h)= Q_0(h),\ \ then \ 
\wt h=h,
\]
\[
\lb{2.27}
l_m(\wt h)\ge l_m(h).
\]

 We show the possibility of this principle  in the
following Lemma.

\begin{lemma} \lb{L2.5}
For each $h\in\ell^\iy$ the estimate \er{2.29} and the following estimate hold:
\[
\lb{2.28}
\nu_n\le h_n, \qqq all \qq n\in\Z.
\]
\end{lemma}
\no {\it Proof.} We apply estimate \er{2.23} to $h$ and to the new sequence:
 $\wt h_m=h_n$ if $m=n$ and $\wt h_n=0$ if $ m\ne n.$  It  is clear
that  $z(k,\wt h)=\sqrt{(k-u_n)^2+h_n^2}$ (the principal value).
Then
$$
y(k,h)\le\Im(\sqrt{(k-u_n)^2+h_n^2}),\ k\in K_+(h).
$$
Then asymptotics \er{2.21} of the function $z(k,h)$ as
 $k\to u_n+ih_n$ yields \er{2.28}.
In  order  to  prove \er{2.29} we use \er{2.26} since
$Q_0(\wt h)=h_m^2/2$. $\BBox$

We recall estimates from [K4].

\begin{theorem} \lb{T3.3}
 Let $h\in \ell^{\iy}$. Then for any  $n\in\Z$ the following estimate holds:
\[
\lb{3.3}
2h_n^2\le\pi \max\lt\{1,{h_n\/r}\rt\}\iint_{D_n(r)}|z'(k)-1|^2dudv,\qq
D_n(r)=(u_n, u_n+r)\times(-h_n, h_n).
\]
If in addition, $\inf_n ( u_{n+1}-u_n)=u_*>0$,  
then
\[
\lb{3.6}
{\pi \/  4} I_D\le\|h\|^2\le {\pi^2 \/  2}\max
\bigg\{1,\,\frac{I_D^{1/2}}{u_*}\bigg\}\,I_D,
\]
\[
\lb{3.8}
{\|l\| \/  2}\le \|J\|\le \sqrt{2}\|l\|(1+{\sqrt{2} \/  u_*}\|l\|).
\]
\end{theorem}

Introduce the domain $S_r=\{z\in\C :\ |\Re z|<r\}, r>0$.
In order to prove Theorem \ref{T3.5} we need the following result about the simple mapping.

\begin{lemma} \lb{L3.4}
The function $f(k)=\sqrt{k^2+h^2},\ k\in\C \sm [-ih,ih],
h>0\,$ is the conformal mapping from $\C\sm [-ih, ih]$ onto
$\C\sm [-h,h]$ and $S_r\sm [-h,h]\subset f(S_r\sm [-ih,ih])$  for any $r>0$.
\end{lemma}
\no {\it Proof.} Consider the image of the half-line $k=r+iv,
v>0$. We have the equations
\[
\lb{3.9}
x^2+y^2=\x \ev r^2+h^2-v^2,\ \ \ \ \  \ \ xy= rv.
\]
The second identity in \er{3.9} yields $x>0$ since $y>0$ . Then
$x^4-\x x^2-r^2v^2=0,$
and enough to check the following inequality
$x^2={1\/2}(\x +\sqrt{\x ^2+4r^2v^2})  > r^2.$
The last estimate follows from the simple relations
$$
(r^2+h^2-v^2)^2+4r^2v^2>(r^2+v^2-h^2)^2,\ \ \ \ \ \ 4r^2v^2>4r^2(v^2-h^2).
\ \ \ \ \BBox
$$

We prove the local estimates for the small slits.
Recall $S_r=\{z\in\C :\ |\Re z|<r\}, r>0$.

\begin{theorem} \lb{T3.5}
Let $h\in \ell^{\iy}$. Assume that
$(u_n-r, u_n+r)\subset(u_{n-1}, u_{n+1})$ and $h_n\le {r\/2}$,
for some $n\in\Z $ and $r>0$. Then
\[
\lb{3.10}
|h_n-|\mu_n^{\pm}||\le \frac{2+\pi}{r}|\mu_n^{\pm}|\sqrt{I_n},
\ \ \ I_n=\frac1{\pi}\int\!\!\!\!\int_{u_n+S_{r}}|z'(k)-1|^2dudv,
\]
\[
\lb{3.11}
0\le h_n-\nu_n\le 2\frac{2+\pi}{r}h_n \sqrt{I_n},
\]
\[
\lb{3.12}
0\le h_n-\frac{l_n}2\le \frac{2+\pi}{r}h_n\sqrt{I_n}.
\]
\end{theorem}
\no {\it Proof.}
Define the functions $f(k)=\sqrt{k^2+h_n^2}$, $\ k\in S_{r}\sm [-ih_n,ih_n],g=f^{-1}$ and $F(w)=z(u_n+g(w),h), w=p+iq$, where the variable
$w\in G_1=f\big((S_{r}\sm [-ih_n,ih_n]\big))$.
The function $F$ is real for real $w$, then $F$ is analytic in the domain
$G=G_1\cup [-h_n,h_n]$ and Lemma \ref{L3.4} yields $S_{r}\ss G$.
Let now $|w_1|={r\/2}$ and $B_r=\{z:|z|<r\}$. Then the following estimates hold
\[
\lb{3.13}
 \sqrt{\pi}\frac{r}2|F'(w_1)-1| \le (\iint_{B_r}
\big|F'(w)-1\big|^2dpdq)^{1/2}\le
\]
$$
 \le\bigg(\iint_{B_r}
\big|(F(w)-g(w))'\big|^2dpdq\bigg)^{1/2}+\bigg(\iint_{B_r}
\big|g'(w)-1\big|^2dpdq\bigg)^{1/2}.
$$
The invariance of the Dirichlet integral with respect to
the conformal mapping gives
\[
\lb{3.14}
 \iint_{B_r}|(F(w)-g(w))'|^2\,dpdq=
\iint_{g(B_r)}|z'(k)-1|^2dudv\le
\iint_{S_{r}+u_n}|z'(k)-1|^2dudv=\pi I_n.
\]
Moreover, the identity $2Q_0=I_D$ implies
\[
\lb{3.15}
{1\/\pi}\iint_{B_r}|g'(w)-1|^2dpdq\le \frac1{\pi}\iint_{\C }|g'(w)-1|^2dpdq
=\frac2{\pi}\int_{-h_n}^{h_n}\sqrt{h_n^2-x^2}\ dx=h_n^2.
\]
Then \er{3.13}-\er{3.15} for $|w_1|={r\/2}$  yields
$
\big|F'(w_1)-1\big|\le\frac2{r}(\sqrt{I_n}+h_n),
$
and \er{3.3} gives
$$
h_n^2\le\frac{\pi^2}4\bigg(\frac1{\pi}\iint_{S_{r}+u_n}|z'(k)-1|^2dudv
\bigg)\le \frac{\pi^2}4 I_n.
$$
Then for $|w_1|=r/2$  we have
\[
\lb{3.16}
|F'(w_1)-1|\le\frac2{r}(1+\frac{\pi}2)\sqrt{I_n}={2+\pi \/  r}\sqrt{I_n},
\]
and the maximum principle yields the needed estimates for $|w_1|\le r/2.$

We prove \er{3.10} for $\mu_n^{+}$. The definition of $\m_n^{\pm}$
(see \er{2.20}) implies
$$
F'(h_n)=\lim_{x\searrow h_n}z'(u_n+g(x))\cdot
g'(x)=\lim_{x\searrow h_n}\frac{g(x)}{\mu_n^{+}}\cdot
\frac x{g(x)}=\frac{h_n}{\mu_n^{+}}.
$$
The substitution of the last identity into \er{3.16} gives \er{3.10}.
The proof for $\m_n^-$ is similar.

We show \er{3.11}.  The definition of  $\nu_n$ (see \er{2.21}) yields
$$
(z(k)-z_n)^2=2i\nu_n(k-u_n-ih_n)(1+o(1)) \ \ as \  \ k\to u_n+ih_n,
$$
$$
g(w)-ih_n=-\frac i{2h_n}(w-z_n)^2(1+o(1))\ \ as \ w\to u_n.
$$
Then we have $F'(0)=\sqrt{\frac{\nu_n}{h_n}}$ and
the substitution of the last identity into \er{3.16} shows
$\bigg|\sqrt{\frac{\nu_n}{h_n}}-1\bigg|\le {2+\pi\/r} \sqrt{I_n}$, which gives \er{3.11}, since by \er{2.30}, $\nu_n\le h_n$. 
Estimate \er{3.16} yiels
$$
0\le2 h_n-l_n=\int_{-h_n}^{h_n}(1-F'(x))dx\le\frac{2h_n}r(2+\pi)\sqrt{I_n},
\qqq  
$$
which implies \er{3.12}. \BBox

\no We prove the estimates in terms of the norm of the space $\ell^p, p\ge 1$.

\no {\bf  Proof of Theorem  \ref{T3.6}.} Estimate \er{1.5} and the H\"older inequality yield \er{3.17}. Using \er{2.29}, $\pi Q_0\le\sum h_nl_n$ (see \er{1.5})
and the H\"older inequality  we obtain
$$
 (\pi /2) I_D\le\|h\|_{\iy}^{1-{p \/  q}}\sum_{n\in\Z }l_n
\,h_n^{p/q}\le(I_D)^{(1-{p \/  q})/2}\|l\|_p\|h\|_p^{{p \/  q}},
$$
$$
 (\pi /2) I_D^{\frac{(1+{p\/q})}2}\le \|l\|_p\|h\|_p^{{p\/q}},
\ \ \ \ \ {\rm and}\ \ \ \ \ \
I_D\le \rt({2\/\pi}\rt)^{{2\/p}}\,\|l\|_p^{{2\/p}}
\|h\|_p^{{2\/q}}.
$$
Estimate \er{3.17} at $q=1$ implies the first one in \er{3.19}.
The last result and \er{2.29} yield the first  inequality in
\er{3.20} and then the second one in \er{3.19}. The second estimate in
\er{3.20} follows from $l_n\le2 h_n, \ n\in\Z $ (see \er{1.3}). 
We have proved \er{2.29} in Lemma \ref{L2.5}.
$\BBox$

\section{Proof of the mains theorems}
\setcounter{equation}{0}


\no {\bf Proof of Theorem \ref{T2.1}}.
Let $1\le p\le 2$ and $r={u_*\/2}$. Estimate \er{3.3} implies
\[
\lb{3.21}
h_n^2\le{\pi^2\/2} \max\{1,{h_n \/  r}\}I_n,  \ \ \
I_n={1\/\pi} \iint_{D_n}\!\!\!|z'(k)-1|^2dudv,\ \ \
D_n=\{ k: |\Re (k-u_n)|<{u_*\/2}\}.
\]
Hence
\[
\lb{3.22}
  h_n\le \frac{\pi^2}{u_*}I_n,\ \ \ if \ \  h_n>{u_* \/  4},
\ \ and \ \
  h_n\le \frac{\pi }{u_*}\sqrt{I_n},\ \ \ if \ \  h_n\le {u_* \/  4}.
\]
Moreover, \er{3.12} yields
\[
\lb{3.23}
h_n\le {l_n\/2}+2h_n {2+\pi\/u_*}\sqrt{I_n},\qqq if  \qq  h_n<{u_*\/4},
\]
and then
\[
\lb{3.24}
h_n\le 2\pi \frac{2+\pi }{u_*}I_n, \qqq if\qq h_n<{u_*\/4}, \ l_n\le h_n,
\]
since $h_n\le {\pi \/2}\sqrt{I_n}$. Hence
$$
if\ \  \ l_n\le h_n, \ \  then \ \
  h_n\le 2\pi \frac{2+\pi}{u_*}I_n=C_1I_n,\ \ \ \ C_1=2\pi {2+\pi\/u_*}.
$$
The last inequality and \er{2.29} yield
\[
\lb{3.25}
\| h\|_p\le
(\sum _{h_n<l_n}h_n^{p})^{1/p}+ (\sum _{l_n\le h_n}h_nh_+^{p-1})^{1/p} \le
\|l\|_p+ C_1^{\frac1{p}}I_D^{\frac{p+1}{2p}}.
\]
If we assume that $C_1^{\frac1{p}}I_D^{\frac{p+1}{2p}}\le \|l\|_p$,
then we obtain $\|h\|_p\le 2\|l\|_p.$

Conversely,  let $\|l\|_p\le C_1^{\frac1{p}}I_D^{\frac{p+1}{2p}}.$
Then \er{3.25}, \er{3.18} implies
$$
\|h\|_p\le
2C_1^{\frac1{p}}
\lt[(2/\pi )^{2/p} \|h\|_{p}^{2/q}\,\|l\|_p^{2/p}\rt]^{\frac{p+1}{2p}}.
$$
Recall that  the  constant $\a_p=(2^{p+2}(2+\pi )/u_*)^p/\pi $.
Hence
$$
\|h\|_p^{1/p^2}\le 2C_1^{1\/  p}\rt[(2 /  \pi)\|l\|_p\rt]^{p+1 \/  p^2},
\ \ \ \ \ {\rm and} \ \ \ \ \ \
\|h\|_p\le 2^{p^2}C_1^p (2/\pi) ^{p+1} \|l\|_p^{1+p},
$$
which yields \er{2.2}. Let $p\ge 2$. Using inequality \er{3.5}, \er{2.29} we obtain
\[
\lb{3.26}
\|h\|_p\le (\sum h_+^{p-2} h_n^2)^{1/p}\le
C_pb^{1/p}I_D^{1/2}, \ \ \ \ b=b(h_+),\  \ h_+=\|h\|_\iy.
\]
Consider the case $b\le 1.$ Then \er{3.26}, \er{3.17} imply
$$
\|h\|_p^2\le C_p^2I_D\le C_p^2(2/\pi) \|h\|_p\|l\|_q,
\ \ \ \ \ {\rm and} \ \ \ \ \ \ \ \|h\|_p\le (2C_p^2/\pi) \|l\|_q,
$$
Consider the case $b>1.$ Then the substitution of \er{2.29}, \er{3.17}
into \er{3.26} yield
$$
\|h\|_p\le C_pI_D^{\frac{p+1}{2p}}u_*^{-1/p}\le
u_*^{-1/p}C_p[(2/\pi) \|h\|_p\|l\|_q]^{\frac{p+1}{2p}},
$$
and
$$
\|h\|_p^{\frac{p-1}{2p}}\le C_pu_*^{-1/p}[(2/\pi) \|l\|_q]^{\frac{p+1}{2p}},
\ \ \ \ \ {\rm and} \ \ \  \ \
\|h\|_p\le (C_pu_*^{-1/p})^{\frac{2p}{p-1}}  (2/\pi)^{\frac{p+1}{p-1}}
 \|l\|_q]^{\frac{p+1}{p-1}},
$$
and combining these two cases we have \er{2.3}.
Inequality $l_n\le 2J_n$ (see \er{2.30}) yields the first estimate in \er{2.4}.
Relation \er{2.30} implies
$$
\|J\|_p^p=\sum |J_n|^p\le \sum (2/\pi )^{p/2}h_n^{p/2} l_n^{p/2}\le
(2/ \pi )^{p/2} \|h\|_p^{p/2} \|l\|_p^{p/2},
$$
and using \er{2.2} we obtain the second estimate in \er{2.4}:
$$
\|J\|_p\le \sqrt{2/\pi }\|l\|_p^{1/2} [2\|l\|_p(1+\a_p\|l\|_p^p)]^{1/2}=
{2\/  \sqrt{\pi }}\|l\|_p(1+\a_p\|l\|_p^p)^{1/2}.
$$
Inequality  $J_n^2\le 4h_n^2/\pi $ (see \er{2.30})    yields
the first estimate in \er{2.5}.
Using \er{2.2} and $\|l\|_p\le 2\|J\|_p$ (see \er{2.4})  we deduce that
$$
\|h\|_p\le 2\|l\|_p(1+\a_p\|l\|_p^p)\le
4\|J\|_p(1+\a_p 2^{p}\|J\|_p^p). \ \ \ \BBox
$$
Consider now the examples, which proves Proposition \ref{T-E},
i.e., the exactness of some estimates.

\no {\bf Example 1.}
Define the sequence $h_n=N,|n|\le N,$ and
$h_n=0, $ if $ |n|>N$ and assume that $u_n=n$,
$\ n\in\Z $. Let $B_r=\{z:|z|<r\}, r>0$. Introduce the function
$g(k)=k+{R^2\/k},\ |k|>R=N\sqrt 8$, which is the conformal mapping
from $\C \bs{\ol {B_R}}$ onto $\C \bs[-2R,2R]$.
Note that $\cup_{|n|\le N}\big[-ih_n+u_n,u_n+ih_n\big]
\ss {\ol{B_R}}$. Then Theorem \ref{T3.3} yields
$\|h\|_{\iy}=N$, $\ u_*=1$, $\ \|h\|_2^2=2N^3$,
$\ \|h\|_1=2N^2$ and using \er{2.29}, \er{2.26} we obtain
$$
N^2=\|h\|_{\iy}^2\le2Q_0=I_D\le16N^2.
$$
Hence inequality \er{3.6} is precise. Moreover, estimate \er{1.5} yields
$\frac{\pi N}2\le\frac{\pi Q_0}N $ and $\|l\|_1\le\frac{2\pi Q_0}N\le16\pi N.$
Then we deduce that  \er{2.2} at $p=1,$ is precise. $\BBox$

In order to consider estimates \er{3.7} we  need the following simple
result.

\begin{lemma} \lb{L3.7}
Let $h\in l^{\iy}, \ u_{1}-u_0=u_0-u_{-1}=1$. Then
\[
\lb{3.27}
l_0\le\big((M^2-h_0^2+1)^2+4h_0^2\big)^{1/4},\ \ \ \
M=\min\big\{h_{-1},h_{1}\big\}.
\]
\end{lemma}
\no {\it Proof.}
Define the sequences
$$
\wt h_m=\cases {
0, &if\ \ $|m|\ge2,$\cr
M, & if\ \ $|m|=1,$\cr
h_n, & if\ \ $ m=0.$\cr }\ \ \ \
\ \ \ \ \ \
\e _m=\cases {
\wt h_m, &if \ \ $m\neq 0,$\cr
      0, &if \ \ \ $m=0.$\cr}
$$
Inequality \er{2.27} implies $\wt l_0\ge l_0$, then enough to show
estimate \er{3.27} only for the sequence $h=\wt h$. We have
$$
z(k,h)=\sqrt{\big(z(k,\e )\big)^2+\big|z(ih_0,\e )\big|^2}.
$$
Hence the maximum principle yields
$$
 l_0\le \Im  z\big(ih_0,\e \big)\le |\sqrt{M^2+(1+ih_0)^2}|\le
 \big((M^2-h_0^2+1)^2+4h_0^2\big)^{1/4}.\ \ \ \ \BBox
$$

\no {\bf Example 2.}
Introduce now the sequences
$$
u_n=n,\ \ n\in\Z,\ \ \ \
h_n=\cases {
N-|n|, &if \  \ $0\le|n|\le N,$\cr
0, &if\ \ $ |n|>N.$\cr}
$$
We estimate  $\|l\|_2$. Using Lemma \ref{L3.7} we obtain for
$|m|\le N-2,\  \x=N-|m|$:
$$
l_m^4\le ( (\x-1)^2-\x^2+1)^2+4\x^2=4(\x-1)^2+4\x^2\le 8N^2.
$$
Moreover, the simple estimate $l_n\le 2h_n$ implies
$l_{N-1}\le 1, \ l_{1-N}\le  1$.
Consider now $A=\sum_{|n|\le N}l_n $. By the Theorem Ivanov-Pomerenke,
 $\frac A4$ is equal the capacity of the compact set
 $ E=\cup_{|n|\le N}[u_n-ih_n,u_n+ih_n]$.
The capacity of the set $ E$ is less than the diameter
which is equal to $2N$ . Then $A\le8N$ and we have
$$
\|l\|^2=\sum l_m^2\le\root{4}\of{8}\,\sqrt N
\sum l_m\le8\root{4}\of{8}\,N^{3/2}=BN^{3/2}.
$$
Assume that the estimate $\|h\|\le C\|l\|(1+\|l\|^p)$ holds
for some constants $C,p>0$. Then for our Example 2  for large
$N$ we obtain  $\|h\|\le 2CB^{p/2}N^{3p/4}$.
On the other hand we have
$\|h\|\ge C_1N^{3/2}$ for some constant $C_1>0$.  Then
 $N^{\frac32(1-p/2)}\le\frac{2C}{C_1}\,B^{p/2}$. It is possible
only for $p\ge2$. Hence estimate $\|h\|\le C\|l\|(1+\|l\|^p)$ is true only for some $p\ge2$ (recall that $p=3$ in \er{3.7}).
$\BBox$

We have considered the estimates for the weight $\o_n\ge 1$.
Now we obtain the counterexample, which shows impossibility of double-sided
estimates in the space $\ell^ \iy _{\o}$ with $\o_n\le 1$.

\no {\bf Counterexample 3.}
Consider now the uniform comb with $u_n=\pi n$,
$\ H=h_n=\|h\|_{\iy}$, $\ n\in\Z $. It is clear
(see [LS]), that in this case $l=l_n=2\arcsin \hbox{th}\,H$,
$\ n\in\Z $. Then $l\to\pi$ as $H\to\iy$ and in this case
for any sequence $\o=(\o_n)_{n\in\Z }$,
$\ \o_n>0$, $\ n\in\Z $; $\ \sum\limits_{n\in\Z }\o_n<+\iy$ and
any $1\le p<+\iy$ the sequence $h$ belongs to $\ell^p_{\o}$ and
$$
\|h\|_{p,\o}^p=\sum_{n\in\Z }h_n^p \o_n=H^p\sum_{n\in\Z } \o_n;
\ \ \ \ \ \ {\rm and }\ \  \ \ \ \
 \|l\|_{p,\o}^p=l^p\sum_{n\in\Z }\o_n.
$$
Hence the estimate  $\|h\|_{p,\o}^{}\le F\big( \|l\|_{p,\o}\big)$
is impossible for some function $F$. Moreover, by the same reason,
the following estimate $\|h\|_{\iy}\le F\big(\|l\|_{\iy}\big)$
 is impossible too. $\BBox$

\no {\bf Proof of  Theorem \ref{T2.4}.}
Recall that for $h\in \ell^\iy, h_n\to 0$ as $n\to\iy$, we define the new sequence $\wt h=\wt h(h)$ by:

\no we take $n_1\in \Z$ such that $\wt h_{n_1}=h_{n_1}=\max_{n\in\Z }h_n>0$; assume that the numbers
$n_1,n_2,\dots,n_k $ have been defined, then we take $n_{k+1}$ such that
$$
\wt h_{n_{k+1}}=h_{n_{k+1}}=\max_{n\in B}h_n>0,\ \ 
 B=\{n\in\Z : |u_n-u_{n_l}|>h_{n_l}, 1\le l\le k\}.
$$
Moreover, we let $\wt h_n=0 $,
if the number $n \notin \{n_k, k\in \Z\}$.
The Lindel\"of principal yields $Q_0(\wt h)\le Q_0(h) $. On the other hand
 open squares $P_k=(u_{n_k}-t_k,u_{n_k}+t_k)\times(-t_k,t_k),
t_k\ev h_{n_k}=\wt h_{n_k}, k\in Z $, does not
overlap. Then applying \er{3.3} to the function
 $(z(k,\wt h)-k) $ and $P_k$, we obtain
\[
\lb{3.29}
2t_k^2\le\pi\int \int_{P_k}|z'(k,\wt h)-1|^2\,dudv,
\]
and
$$
Q_0(\wt h)=\frac12I_D(\wt h)\ge\frac1{\pi^2}\sum_{k\ge1}t_k^2=\frac1{\pi^2}
\|\wt h\|^2_2
$$
which yields the first estimate in \er{2.15}.

   Let $\O_k=\{n\in\Z : u_n\in[u_{n_k}-t_k,u_{n_k}+t_k]\} $.
By the Lindel\"of principal, the  gap length $l_{n_k}$ such that
 $[u_n,u_n+ih_n], n\in \O_k $ increases if we take off all another slits.
By the Theorem of Ivanov-Pomerenke (see Section 1), the sum of new gap lengths equals to $4\times capacity $ of the set
 $E=\cup_{n\in \O_k}[u_n-h_n,u_n+h_n] $, which is less than
the diameter of the set $E$. Then
$\sum_{n\in \O_k}l_n(h)\le2\sqrt2t_k$,
and using the last estimate we obtain
\[
\lb{3.31}
\pi Q_0(h)\le\sum_{n\in\Z }h_nl_n\le \sum_{k\ge1}\sum_{n\in \O_k}h_nl_n \le
\sum_{k\ge1}t_k\sum_{n\in \O_k}l_n\le 2\sqrt2\sum_{k\ge1}t_k^2=2\sqrt2\|\wt h\|^2_2,
\]
since $h_n\le t_k, n\in \O_k $ and the diameter of the set $E $
is less than or equals $2\sqrt2t_k $.  $\BBox$

Note that the proved Theorem shows that estimates \er{3.7},\er{3.8}
are fulfilled for the weaker conditions on the sequence $u_n,n\in \Z$.

Recall the following identity for $v(z)=\Im k(z), z=x+iy$ from [KK1]:
\[
\lb{3.32}
v(x)=v_n(x)\big(1+V_n(x)\big), \quad
V_n(x)=\frac1{\pi}\!\!\! \intl_{\R \bs \g_n}\frac{v(t)dt}{|t-x|\,v_n(t)}, \quad
v_n(x)=|(x-z_n^{+})(x-z_n^{-})|^{1\/2},
\]
for all $x\in \g_n=(z_n^{-},z_n^{+})$. 
In order to prove Theorem \ref{T2.2} we need the following results. 

\begin{lemma} \lb{L3.8}
Let $h\in \ell^\iy $ and $u_*>0$ and $\x=e^{\|h\|_{\iy}\/u_*}$. Then
the following estimates hold:
\[
\lb{3.33}
s\ev s(h)\ev \inf |\s_n(h)|\le u_* \le
{\pi s\/2}\max\lt\{e^2,\x^{5\pi \/2}\lt\},
\]
\[
\lb{3.34}
1+{2\|h\|_{\iy}\/  s\pi}\le \x^9 ,
\]
\[
\lb{3.35}
\max_{n\in \g_{n} } V_n(x)\le {2\|h\|_{\iy}\/  \pi s}, \ \ n\in \Z,
\]
\[
\lb{3.36}
2h_n\le l_n(1+\max_{n\in \g_{n} } V_n(x))\le l_n(1+{2\|h\|_{\iy} \/  \pi s})
\le l_n \x^9, \ \ n\in \Z.
\]
\end{lemma}
\no {\it Proof.} Introduce the domain
$$
G=\big\{z\in\C : h_+\ge\Im\,z>0,\ \Re\,z
\in(-{u_* \/  2},{u_* \/  2})\big\}\cup\{\Im\,z>h_+\big\}, \ \
h_+=\|h\|_{\iy}.
$$
Let $g$ be the conformal mapping from $G$ onto $\C _+$, such that
$g(iy)\sim iy$ as $y\nearrow +\iy$
and let $\a , \b $ be images of the points ${u_*\/2}, {u_*\/2}+ih_+$
respectively. Define the function $f=\Im g$. Fix any $n\in \Z $.
Then the maximum principle
yields
$$
y(k)=\hbox{Im}\,z(k,h)\ge f(k-p_n),\quad k\in G+p_n,\ \ \
p_n={1 \/  2}(u_{n-1}+u_n).
$$
Due to the fact that these positive functions equal zero
on the interval $(p_n- {u_*\/2},p_n+{u_*\/2})$, we obtain
$$
\frac{\pa y}{\pa v}\,(x)=\frac{\pa x}{\pa u}\,(x)\ge
\frac{\pa f}{\pa v}\,(x-p_n),\quad x\in \big(p_n-\frac{u_*}2,
p_n+\frac{u_*}2\big).
$$
Then
$$
z(u_{n})-z(u_{n-1})\ge\intl_{-u_*/2}^{u_*/2}\frac{\pa f}{\pa v} \,(x)\,dx=
2\a >0,
$$
and the estimate $s\le u_*$ (see [KK1]) implies $2\a\le s\le u_*$.
Let $w: \C_+\to G$ be the inverse function for $g$, which is
defined uniquely and the Christoffel-Schwartz formula  yields
$$
w(z)=\int _0^z\sqrt{{t^2-\b^2 \/  t^2-\a^2}}dt,\ \ \ \ \ 0<\a<\b.
$$
Then we have
\[
\lb{3.37}
{u_* \/  2}=\int _0^{\a}\sqrt{{\b^2 -t^2\/  \a^2-t^2}}dt,
\ \ \
h_+=\int _{\a}^{\b}\sqrt{{\b^2-t^2 \/  t^2-\a^2}}dt.
\]
The first integral in \er{3.37} has the simple double-sided estimates
$$
\a=\int _{0}^{\a}dt\le {u_* \/  2}\le
\int _0^{\a}{\b dt\/  \sqrt{\a^2-t^2}}={\b\pi \/  2},
$$
that is
\[
\lb{3.38}
2\a\le s\le u_*\le \pi \b.
\]
Consider the second integral in \er{3.37}. Let $\ve =\b /\a\ge 5$ and using
the new variable $t=\a \cosh r, \cosh \d =\ve $, we obtain
$$
h_+=\a \int _0^{\d }\sqrt{\ve ^2-\cosh ^2r}dr\ge
\a\ve\int_0^{\d /2}\sqrt{1-{\cosh ^2r \/  \ve ^2}}dr\ge \b\d{2\/  5},
$$
since for $r\le \d/2$ we have the simple inequality
$$
{\cosh ^2r \/  \cosh ^2\d}  \le e^{-\d} (1+e^{-\d})^2 \le \ve^{-1}
(1+\ve^{-1})^2.
$$
Due to $\ve \le e^{\d}$ we get $\ve \le \exp (5h_+/2\b)$
and estimate \er{3.38} implies
\[
\lb{3.39}
{1 \/  s}\le {\pi  \/  2u_*}\exp ({5\pi  \/  2u_*}h_+),
\ \ \ if  \ \ \ \ve \ge 5.
\]
If $\ve \le 5$, then using \er{3.38} again we obtain
$$
{1 \/  s}\le {\ve  \/  2\b}\le {\pi \ve  \/  2u_*}\le
{\pi \/  2u_*}5, \ \ \ \ \ \    if  \ \ \ \ve \le 5.
$$
and the last estimate together with \er{3.39} yield \er{3.33},\er{3.34}.

 Identity \er{3.32} for $ x\in \g_n=(z_n^{-},z_n^{+})$ implies
$$
\pi V_n(x)=
\int_{-\iy}^{z_n^--s}\frac{v(t)dt}{|t-x|v_n(t)} +
\int_{z_n^++s}^{\iy}\frac{v(t)dt}{|t-x|v_n(t)}\le
\int_{-\iy}^{z_n^--s}\frac{h_+dt}{|t-z_n^-|^2} +
\int_{z_n^++s}^{\iy}\frac{h_+dt}{|t-z_n^+|^2}\le
 {2h_+\/s}.
$$
Using \er{3.32}, \er{3.34} and simple inequality
$v_n(z_n)\le l_n/2$ we have \er{3.36}.  \ \ $\BBox$

We prove the two-sided estimates  of $h_n, l_n, \m_n^{\pm}, J_n$ in
the weight spaces.

\no {\bf  Proof of Theorem \ref{T2.2}}.
 The first estimate in \er{2.6} follows from
$h_n\le 2\pi |\m_{n}^{\pm}|$ (see [KK1]). The second one in \er{2.6} follows from
 $\|h\|_{\iy}\le \sqrt{I_D}=\|J\|\le \|J\|_p\le \|J\|_{p,\o}$
since $\o_n\ge 1$ for any $n\in \Z$.
Moreover, substituting \er{3.18} into $\|h\|_{\iy}\le \sqrt{I_D}$,
using \er{2.2} and $\|f\|_p\le \|f\|_{p,\o}$  for any $f$,  we obtain  the last estimate in \er{2.6}.

Introduce the function $\x=\exp {{\|h\|_{\iy }\/u_*}}$.
The first estimate in \er{2.7} follows from \er{1.3}.
Due to \er{3.36} we  get $2h_n\le \x^9l_n$, which yields the second
estimate in \er{2.7}.

The first estimate in \er{2.8} follows from \er{2.30}.
Using \er{2.30}, \er{3.36} we  have $J_n^2\le 2l_nh_n/\pi \le (\x^9/\pi )l_n^2$, which gives the second estimate in \er{2.8}.

The first estimate in \er{2.9} follows from \er{2.30}, \er{1.3}.
Using \er{2.30}, \er{3.36} we obtain $h_n^2\le \x^9l_nh_n/2 \le (\pi \x^9/2)J_n^2$, which yields the second inequality in  \er{2.9}.

Identity $2\m_n^{\pm}=\pm l_n[1+V_n(z_n^{\pm})]^2$  (see [KK1])
implies $2|\m_n^{\pm}|\ge l_n$, which yields
the first inequality in  \er{2.10}. Moreover, using \er{3.36} we obtain
the estimate $2|\m_n^{\pm}|\le \x^{18}l_n$, which  gives
the second one in \er{2.10}.\ \ $\BBox$

\no {\bf Acknowledgments.}
E. Korotyaev was partly supported by DFG project BR691/23-1.
The various parts of this paper were written at the Mittag-Leffler Institute, Stockholm  and in the Erwin Schr\"odinger Institute for Mathematical Physics, Vienna, E. Korotyaev is grateful to the Institutes for the hospitality.

\end{document}